\documentclass[12pt,english]{amsart}
\usepackage[T1]{fontenc}
\usepackage[latin9]{inputenc}
\usepackage{geometry}
\usepackage{babel}
\usepackage{verbatim}
\usepackage{amsthm}
\usepackage{amstext}
\usepackage{amssymb}
\usepackage[unicode=true]
 {hyperref}
\usepackage{breakurl}

\makeatletter
\numberwithin{equation}{section}
\numberwithin{figure}{section}
\theoremstyle{plain}
\newtheorem{thm}{\protect\theoremname}[section]
  \theoremstyle{remark}
  \newtheorem{rem}[thm]{\protect\remarkname}

\@ifundefined{definecolor}
 {\usepackage{color}}{}

\@ifundefined{definecolor}{\usepackage{color}}{}
\usepackage{amsfonts}\usepackage{babel}

\providecommand{\U}[1]{\protect\rule{.1in}{.1in}}

\usepackage{calrsfs}
\DeclareMathAlphabet{\pazocal}{OMS}{zplm}{m}{n}

\usepackage{babel}

\makeatother

  \providecommand{\remarkname}{Remark}
\providecommand{\theoremname}{Theorem}

\begin{document}

\title[Interpolation of compact operators]{Nigel Kalton and complex interpolation of compact operators}

\author{{Michael Cwikel}}

\address{Cwikel: Department of Mathematics, Technion - Israel Institute of
Technology, Haifa 32000, Israel }

\email{mcwikel@math.technion.ac.il }

\author{{Richard Rochberg}}

\address{Rochberg: Department of Mathematics, Washington University, St. Louis,
MO, 63130, USA }

\email{rr@math.wustl.edu}

\thanks{The first named author's work was supported by the Technion V.P.R.\ Fund
and by the Fund for Promotion of Research at the Technion. The second
named author's work was supported by the National Science Foundation
under Grant No. 1001488.}
\begin{abstract}
The fundamental theorem of complex interpolation is 

\textbf{Boundedness Theorem:} If, for $j=0,1$, a linear operator
$T$ is a bounded map from the Banach space $X_{j}$ to the Banach
space $Y_{j}$ then, for each $\theta$, $0<\theta<1$, $T$ is a
bounded map between the complex interpolation spaces $\left[X_{0},X_{1}\right]_{\theta}$
and $\left[Y_{0},Y_{1}\right]_{\theta}$. 

Alberto Calder\'on, in his foundational presentation of this material
fifty-one years ago \cite{CalderonA1964}, also proved the following
companion result: 

\textbf{Compactness Theorem/Question:} Furthermore \textbf{\textit{in
some cases}}, if $T$ is also a compact map from $X_{0}$ to $Y_{0}$,
then, for each $\theta$, $T$ is a compact map from $\left[X_{0},X_{1}\right]_{\theta}$
to $\left[Y_{0},Y_{1}\right]_{\theta}$. 

The fundamental question of exactly which cases could be covered by
such a result was not resolved then, and is still not resolved. 

The paper \cite{aaaCwKa}, which is the focus of this commentary,
is a contribution to that question. We will not summarize in any detail
here the contents of \cite{aaaCwKa} or of related works. (Some of
that may be done later in \cite{CwikelMRochbergR2014Preprint}.) Rather
we will take this opportunity to sketch the mathematical world, historical
and current, in which that paper lives. We will see that there have
been many very talented contributors and many fine contributions;
however the core problem remains open. We will at least be able to
announce some small new partial results. 
\end{abstract}
\maketitle
This paper is the fourth (and, together with its technical sequel
\cite{CwikelMRochbergR2014Preprint}, probably the last) in a series
\cite{CwikelMMilmanMRochbergR2014,CwikelMMilmanMRochbergR2014B,CwikelMMilmanMRochbergR2014C}
which we have recently written (the earlier ones in collaboration
with Mario Milman) to survey and discuss some small parts of the very
extensive and very impressive body of research created by our brilliant
colleague and dear friend Nigel Kalton. We have also submitted material
from these papers to possibly appear in a forthcoming ``Selecta''
volume to be published in memory of Nigel%
\footnote{In this paper, as in the earlier papers of this series, we allow ourselves
the informality of addressing Nigel Kalton by his first name. This
should be understood as part of our expression of our great admiration
of his work, of our warm friendship for him, and our deep regret that
he is no longer with us. %
}. 

The topic which we are about to discuss here has an extensive and
varied bibliography. We will certainly not be able to do full justice
to the contributions of all who have worked on it and on other closely
related topics. We are quite possibly not aware of some of these contributions.
We apologize in advance for all omissions. 

We are going to be dealing with lots of different Banach spaces, including
$L^{p}$ spaces. Let us specify at the outset that all of them will
be over the \textit{complex} field. Thus, also, the terminology Banach
couple (see below), will always refer here to a couple of \textit{complex}
Banach spaces. To avoid any ambiguity, we should also mention that,
for any Banach space $X$, we will be using the notation $\ell^{p}(X)$
or $c_{0}(X)$  to denote the Banach space, equipped with its natural
norm, of all two-sided $X$ valued sequences $\left\{ x_{n}\right\} _{n\in\mathbb{Z}}$
for which the two-sided numerical sequence $\left\{ \left\Vert x_{n}\right\Vert \right\} _{n\in\mathbb{Z}}$
is in $\ell^{p}$ or $c_{0}$ respectively. 

We have begun writing a companion paper or set of lecture notes \cite{CwikelMRochbergR2014Preprint},
which is intended to be a sequel to this commentary. We plan to post
a preliminary version of it on the arXiv concurrently with this document.
It will deal with issues which constraints of space and time have
prevented us from including here. It will contain some minor new results,
and discussions of some technical aspects of \cite{aaaCwKa}, including,
for example, some more explicit details for the proof of the result
of that paper dealing with arbitrary couples of Banach lattices.

\section{\label{sec:intro}The history of the problem before the writing of
\cite{aaaCwKa}}

Our story begins with the classical Riesz-Thorin Theorem. According
to that theorem, any linear operator $T:L^{p_{0}}+L^{p_{1}}\to L^{q_{0}}+L^{q_{1}}$
which is a bounded map from $L^{p_{0}}$ into $L^{q_{0}}$ and also
from $L^{p_{1}}$ into $L^{q_{1}}$ must also be a bounded map from
$L^{p}$ into $L^{q}$ whenever the two exponents $p$ and $q$ are
related by the convexity formulae 
\begin{equation}
1/p=\left(1-\theta\right)/p_{0}+\theta/p_{1}\text{ and }1/q=(1-\theta)/q_{0}+\theta/q_{1}\label{eq:convform}
\end{equation}
for some $\theta\in(0,1)$. In 1926 Marcel Riesz gave the first proof
of a somewhat limited version of this theorem. (One could claim that
he was more or less ``forced'' to obtain this theorem after he found
a clever trick to show that the Hilbert transform is bounded on $L^{p}$,
but only when $p$ is an even integer.) In 1939 Olof Thorin published
a quite different new proof of the complete version, using holomorphic
$L^{p}$ valued functions, a device which John Edensor Littlewood
once described as ``the most impudent idea in mathematics''. (Even
those of us who have the greatest admiration for this wonderful proof
might feel obliged to concede that Littlewood was perhaps just a little
too carried away when he wrote that.) 

In 1960 Mark Krasnosel'skii investigated the interplay of compactness
with the Riesz-Thorin Theorem. He showed \cite{KrasnoselskiiM1960}
that if the above-mentioned operator $T$ has the additional property
that it maps $L^{p_{0}}$ compactly into $L^{q_{0}}$, then it also
maps $L$$^{p}$ compactly into $L^{q}$. 

In the mid 1960's Selim Krein, Jacques-Louis Lions and Alberto Calder\'on
independently found ways to upgrade Thorin's ``impudence'' from
the setting of $L^{p}$ spaces into the setting of \textit{Banach
couples}%
\footnote{Some papers use the terminology \textit{interpolation pair} or \textit{Banach
pair}, instead of Banach couple.%
}, i.e., general pairs $\left(X{}_{0},X_{1}\right)$ of Banach spaces
$X_{0}$ and $X_{1}$ which are both continuously embedded in some
Hausdorff topological vector space. 

The theory developed by Alberto Calder\'on had a somewhat wider scope
than those of Krein and Lions. As described in \cite{CalderonA1964}
(and also in many other papers and monographs, including \cite{BerghLofstrom},
\cite{KPS}, \cite{TriebelH1995}, etc.) for each Banach couple $\left(X_{0},X_{1}\right)$,
Calder\'on's construction uses certain holomorphic $X_{0}+X_{1}$
valued functions to obtain, for each $\theta\in(0,1)$, a Banach space,
usually denoted by $\left[X_{0},X_{1}\right]_{\theta}$ and referred
to as a \textit{complex interpolation space}. This space almost automatically
has the following Riesz-Thorin-like \textit{interpolation} property:
Given any two Banach couples $\left(X_{0},X_{1}\right)$ and $\left(Y_{0},Y_{1}\right)$,
any linear operator $T:X_{0}+X_{1}\to Y_{0}+Y{}_{1}$ whose restriction
to $X_{j}$ is a bounded map $T:X_{j}\to Y_{j}$ for $j=0,1$ is also
a bounded map from $\left[X_{0},X_{1}\right]_{\theta}$ to $\left[Y_{0},Y_{1}\right]_{\theta}$
for each $\theta\in(0,1)$. 

The first example of a Banach couple that comes to mind is a pair
of $L^{p}$ spaces, say $\left(L^{p_{0}},L^{p_{1}}\right)$. (Of course
they must both have the same underlying measure space.) Calder\'on
showed that 
\begin{equation}
\left[L^{p_{0}},L^{p_{1}}\right]_{\theta}=L^{p},\label{eq:lpcase}
\end{equation}
where $p$, $p_{0}$, $p_{1}$ and $\theta$ are connected by the
first of the above two convexity formulae (\ref{eq:convform}). We
refer e.g.~to \cite{BerghLofstrom} pp.~106--107 for a quick direct
proof of (\ref{eq:lpcase}) which is essentially a rewriting of Thorin's
proof of the Riesz-Thorin theorem. In Calder\'on's original presentation
of these matters \cite{CalderonA1964}, this result emerged as a particular
example (using Section 13.3 of \cite{CalderonA1964}) of his simple
formula identifying $\left[X_{0},X_{1}\right]_{\theta}$ as the ``product''
$X_{0}^{1-\theta}X_{1}^{\theta}$ (\cite{CalderonA1964} Section 13.6
pp.~124--125) in the case where $\left(X_{0},X_{1}\right)$ is a
\textit{lattice couple}, i.e., when $X_{0}$ and $X_{1}$ are Banach
lattices of complex valued measurable functions on the same measure
space%
\footnote{In fact, for the formula $\left[X_{0},X_{1}\right]_{\theta}=X_{0}^{1-\theta}X_{1}^{\theta}$
to hold, the lattices $X_{0}$ and $X_{1}$ have to satisfy a certain
mild condition which we shall not describe here. In fact Calder\'on
obtains a rather more general ``vector-valued'' version of this
formula where the measurable functions take values in Banach spaces.%
}. Some of the main ideas for the proof of this characterization formula
(\cite{CalderonA1964} Section 33.6 pp.~171--180), also came naturally
from Thorin's proof of the Riesz-Thorin Theorem. 

Although it is of course explicitly defined, the process for obtaining
$\left[X_{0},X_{1}\right]_{\theta}$ from $X_{0}$ and $X_{1}$ when
these two spaces are not $L^{p}$ spaces, or do not have some kind
of compatible lattice structure, works in ways which seem very mysterious.
These ways bewildered even Alberto Calder\'on himself, as he once
``confessed'' to Mitchell Taibleson%
\footnote{In November 1999, immediately after one of us gave a lecture at Washington
University about our forthcoming joint paper \cite{ckmr2002} with
Nigel and Mario Milman, Mitch spontaneously rose to his feet and gave
us all a touching account of a stroll he had taken in Hyde Park, Chicago
with Alberto Calder\'on in the 1960s, just after Alberto had given
a lecture about the exciting new interpolation method that he had
just developed. This was the occasion on which Alberto made the ``confession''
mentioned above. %
}. Some partial compensation for this mystery would be found some decades
later when Svante Janson \cite{JansonS1981}, motivated by some ideas
of Jaak Peetre, revealed, among many other interesting results, that
the complex interpolation method is a special concrete case, involving
Fourier series, of an abstract interpolation method going back to
Aronszajn and Gagliardo. 

Calder\'on's remarkable exposition \cite{CalderonA1964} of his complex
interpolation method contains many ``gems''. They include a beautiful
theorem (in Sections%
\footnote{Here, and also elsewhere, we need to refer to two different sections
of \cite{CalderonA1964} because the format of that paper is usually
that a result is announced in Section $N$ and its proof is given
in Section $N+20$.%
} 12.1, p~121 and 32.1, pp.~148--156) characterizing the dual space
of $\left[X_{0},X_{1}\right]_{\theta}$ (cf.~a slightly different
way of presenting it in \cite{CwikelDualityLectures}) and, in its
deeper depths (in Sections 34.1 and 34.2, pp.~180--188), is the development
and application of the formula still referred to as the ``Calder\'on
Reproducing Formula'' which is one of the forerunners of modern wavelet
analysis, where the formula also survives under the name ``continuous
wavelet transform''. 

But among those ``gems'', the one that mainly concerns us here is
the following generalization (implicit in Sections 9.4 p.~118 and
29.4, pp.~137--138) of Krasnosel'skii's compactness theorem: If the
Banach couple $\left(Y_{0},Y_{1}\right)$ satisfies a certain ``natural''
condition, and if the above-mentioned linear operator $T$ which maps
the space $X_{j}$ to $Y_{j}$ boundedly for $j=0$ and $j=1$ is
also compact for at least one of these values of $j$, say $j=0$,
then $T$ is also compact as a map from $\left[X_{0},X_{1}\right]_{\theta}$
into $\left[Y_{0},Y_{1}\right]_{\theta}$ for each $\theta\in(0,1)$.
The ``natural'' condition which Calder\'on had to impose in this
theorem is, roughly speaking, that $Y_{0}$ can be \textquotedblleft{}approximated\textquotedblright{}
by some sequence or net of its finite dimensional subspaces in a way
which is suitably consistent with $Y_{1}$. 

To streamline the rest of this discussion, we will generally use the
following notation (cf.~\cite[p.~22]{CwikelIntoLattices}): Given
two Banach couples $\left(X_{0},X_{1}\right)$ and $\left(Y_{0},Y_{1}\right)$,
we may write 
\begin{equation}
\left(X_{0},X_{1}\right)\blacktriangleright\left(Y_{0},Y_{1}\right)\label{eq:dbtr}
\end{equation}
to mean that $T:\left[X_{0},X_{1}\right]_{\theta}\to\left[Y_{0},Y_{1}\right]_{\theta}$
is compact for every $\theta\in(0,1)$ whenever $T:X_{0}+X_{1}\to Y_{0}+Y_{1}$
is a linear operator for which $T:X_{0}\to Y_{0}$ compactly and $T:X_{1}\to Y_{1}$
boundedly. Thus, for example, Krasnosel'skii's theorem can be expressed
by writing $\left(L^{p_{0}},L^{p_{1}}\right)\blacktriangleright\left(L^{q_{0}},L^{q_{1}}\right)$. 

We may then also use the notation 
\[
(\ast,\ast)\blacktriangleright\left(Y_{0},Y_{1}\right)
\]
to mean that $\left(Y_{0},Y_{1}\right)$ is a Banach couple with the
property that $\left(X_{0},X_{1}\right)\blacktriangleright\left(Y_{0},Y_{1}\right)$
for \textit{all} Banach couples $\left(X_{0},X_{1}\right)$. Thus,
Calder\'on's result in \cite{CalderonA1964} can be expressed by
stating that $(\ast,\ast)\blacktriangleright\left(Y_{0},Y_{1}\right)$
for every Banach couple $\left(Y_{0},Y_{1}\right)$ which has the
particular approximation property mentioned above. 

Analogously it will sometimes be convenient to write $\left(X_{0},X_{1}\right)\blacktriangleright(\ast,\ast)$
to mean that $\left(X_{0},X_{1}\right)$ is a Banach couple with the
property that $\left(X_{0},X_{1}\right)\blacktriangleright\left(Y_{0},Y_{1}\right)$
for \textit{all} Banach couples $\left(Y_{0},Y_{1}\right)$. Finally,
we can use the notation $(\ast,\ast)\blacktriangleright(\ast,\ast)$
as shorthand for the (thus far unconfirmed) statement that the following
question has an affirmative answer.

\noindent \textbf{\textit{Question CIC}}%
\footnote{CIC$=$Complex interpolation of compactness.%
}\textbf{\textit{: Does the property $\left(X_{0},X_{1}\right)\blacktriangleright\left(Y_{0},Y_{1}\right)$
hold for all Banach couples $\left(X_{0},X_{1}\right)$ and $\left(Y_{0},Y_{1}\right)$?}}

Calder\'on was evidently not able to answer this question. But then
no one else has been able to do this either, in the more than half
a century since Calder\'on submitted his paper \cite{CalderonA1964}
to Studia Mathematica. Quite a number of us have made valiant efforts
to do so, and some of the world's best mathematicians, (including
Alberto Calder\'on himself!%
\footnote{He discussed it briefly with one of us during a visit to the Technion
in 1989.%
}) have been invited to consider (or reconsider!) this problem. Several
of them, including Nigel Kalton, have expressed the opinion that the
answer to Question CIC is ``no''. But where is the counterexample? 
\begin{rem}
\label{rem:TwoSided}One would expect it to be easier to obtain an
affirmative answer to the variant of Question CIC for ``two-sided
compactness'', i.e., where instead of showing that $\left(X_{0},X_{1}\right)\blacktriangleright\left(Y_{0},Y_{1}\right)$
holds, one only seeks to show that a weaker property holds, in which
$T:\left[X_{0},X_{1}\right]_{\theta}\to\left[Y_{0},Y_{1}\right]_{\theta}$
is required to be compact, only for those operators $T$ which satisfy
the additional condition that $T:X_{1}\to Y_{1}$ is also compact.
But this supposedly easier question is also still open. 
\end{rem}
Quite a number of partial results have been obtained during this said
half century. Let us here briefly mention a few%
\footnote{Please keep in mind the apology at the beginning of this document.%
} of them which preceded the particular paper which we have been asked
to discuss. (Some other more recent results will be mentioned later.)

Arne Persson \cite{PerssonA1964} proved a variant of Calder\'on's
result, namely that $(\ast,\ast)\blacktriangleright\left(Y_{0},Y_{1}\right)$
holds also when $\left(Y_{0},Y_{1}\right)$ satisfies a kind of approximation
property which is somewhat different from that imposed by Calder\'on.
A remarkable paper \cite{cobos-peetre} by Fernando Cobos and Jaak
Peetre, even though it did not directly deal with this question, injected
some powerful new ideas into the game. Helped by those ideas and by
others in \cite{JansonS1981}, one of us \cite{CwikelM1992RCI} was
able to show, essentially%
\footnote{This is implicit in Step 1 of the proof of Theorem 2.1 on pp.~339-340
of \cite{CwikelM1992RCI}. %
}, that the problem of answering Question CIC is equivalent to determining
whether or not $\left(X_{0},X_{1}\right)\blacktriangleright\left(Y_{0},Y_{1}\right)$
holds for just one special particular choice of the couples $\left(X_{0},X_{1}\right)$
and $\left(Y_{0},Y_{1}\right)$. 

In order to prepare for a precise formulation of this result, we first
need to recall the definitions of the Banach sequence spaces $FL^{p}$
and their weighted counterparts $FL_{\alpha}^{p}$ (sometimes also
denoted by $FL^{p}(e^{\alpha\nu})$). For each $p\in[1,\infty]$,
the space $FL^{p}$ consists of all sequenences $\left\{ \lambda_{n}\right\} _{n\in\mathbb{Z}}$
of complex numbers which arise as the Fourier coefficients of some
function $f:\mathbb{T}\to\mathbb{C}$ in $L^{p}\left(\mathbb{T}\right)$.
Then, for each $\alpha\in\mathbb{R}$, $FL_{\alpha}^{p}$ consists
of all sequences $\left\{ \lambda_{n}\right\} _{n\in\mathbb{Z}}$
for which $\left\{ e^{\alpha n}\lambda_{n}\right\} _{n\in\mathbb{Z}}$
is an element of $FL^{p}$. Your first guess as to how these spaces
are normed will be correct. These spaces, for $p=1$ and $p=\infty$,
play crucial roles in Svante Janson's important alternative characterisations
\cite{JansonS1981} of complex interpolation spaces, which we mentioned
above. For later purposes we should also mention the subspaces $FC$
and $FC_{\alpha}$ of $FL^{\infty}$ and $FL_{\alpha}^{\infty}$ which
are obtained when $L^{\infty}$ is replaced by its subspace of continuous
functions. 

We can now state the above-mentioned result of \cite{CwikelM1992RCI}
explicitly: 

\textit{The answer to Question CIC is affirmative if and only if the
property} 
\begin{equation}
\left(\ell^{1}\left(FL^{1}\right),\ell^{1}\left(FL_{1}^{1}\right)\right)\blacktriangleright\left(\ell^{\infty}\left(FL^{\infty}\right),\ell^{\infty}\left(FL_{1}^{\infty}\right)\right)\label{eq:feid}
\end{equation}
\textit{holds.} 

\noindent (A more explicit and detailed proof of this equivalence
would appear later in \cite[pp.~355--358]{CwikelMKrugljakNMastyloM1996}.)

Although (\ref{eq:feid}) and Question CIC remained unresolved in
\cite{CwikelM1992RCI}, the interplay between them led to another
result: Without imposing any approximation assumption on $\left(Y_{0},Y_{1}\right)$,
even if compactness perhaps does not ``interpolate'', then at least
it ``extrapolates''. I.e., if $T:X_{j}\to Y_{j}$ is bounded for
$j=0,1$ and if $T:\left[X_{0},X_{1}\right]_{\theta}\to\left[Y_{0},Y_{1}\right]_{\theta}$
is compact for just one value of $\theta$ in $(0,1)$, then it is
compact for all $\theta\in(0,1)$. A different proof of this result
was immediately obtained by Fernando Cobos, Thomas K\"uhn and Tomas
Schonbek \cite{CobosKuhnSchonbek}. They also proved that $\left(X_{0},X_{1}\right)\blacktriangleright\left(Y_{0},Y_{1}\right)$
holds, without the requirement of Calder\'on's or Persson's approximation
hypotheses, in the case where $(X_{0},X_{1})$ and $(Y_{0},Y_{1})$
are both lattice couples (subject to some mild additional conditions). 

Of course, more or less in parallel with the development of the complex
interpolation method, another important method, the so-called real
interpolation method, was also developed in the 1960's, one of the
main steps in that process being the work \cite{LionsJPeetreJ1964}
of Jacques-Louis Lions and Jaak Peetre. It was natural to ask a question
analogous to Question CIC, where, instead of the spaces $\left[X_{0},X_{1}\right]_{\theta}$
and $\left[Y_{0},Y_{1}\right]_{\theta}$ one considers the Lions-Peetre
spaces generated by this method, and which are usually denoted by
$\left(X_{0},X_{1}\right)_{\theta,p}$ and $\left(Y_{0},Y_{1}\right)_{\theta,p}$.
The answer to this question was also not immediately evident. Affirmative
answers for some of its special cases finally led the way to an affirmative
answer for its general form. This can in fact be found in the same
papers \cite{CobosKuhnSchonbek} and \cite{CwikelM1992RCI} which
we mentioned just above in connection with Question CIC. Some sequels
to that answer will be briefly discussed below in Section \ref{sec:RelatedResults}.

\section{The paper \cite{aaaCwKa} itself}

The second named author of this document once good-naturedly claimed
that the first named author is afflicted with a ``disease'' of wrestling
with Question CIC, and from time to time he ``infects'' others with
the same ``disease''. The first named author does not deny this
claim, even if in his older years the ``disease'' has become somewhat
milder. Perhaps his greatest success was that he ``infected'' Nigel
Kalton, at least temporarily, and this led Nigel to make very substantial
contributions in \cite{aaaCwKa}. Of course Nigel's extremely broad
and deep and intensive interests and activities in so many topics
would guarantee him immunity from succumbing to this disease for any
period significantly beyond the time during which he worked with the
first named author on \cite{aaaCwKa}. 

Here would be the natural place in this document for a summary of
the main results and techniques of \cite{aaaCwKa}, including a list
of the various (quite large) classes of Banach couples $\left(X_{0},X_{1}\right)$
and $\left(Y_{0},Y_{1}\right)$, not previously investigated, for
which it is shown there that $\left(X_{0},X_{1}\right)\blacktriangleright\left(Y_{0},Y_{1}\right)$
holds. However the paper speaks quite well for itself in these matters,
and exactly that summary can be seen on page 262 of \cite{aaaCwKa}
(or on page 2 of the preliminary arXiv version).

Let us mention that one of the results of \cite{aaaCwKa}, namely
that 
\begin{equation}
\left(X_{0},X_{1}\right)\blacktriangleright(\ast,\ast)\mbox{ for every lattice couple }\left(X_{0},X_{1}\right),\label{eq:felc}
\end{equation}
is obtained via an interplay (see pages 269--270 of \cite{aaaCwKa})
with a different interpolation method which probably deserves to be
better and more widely known. It is often called the ``plus-minus''
method%
\footnote{The reason for this name comes from a property of (weakly) unconditionally
convergent series which is mentioned, for example, on line 20 of \cite[p.~58]{JansonS1981}.%
}. It comes in two ``flavors''. In its original version this method
was invented by Jaak Peetre \cite{PeetreJ1971} and, in a modified
version which was motivated by applications to Orlicz spaces, it first
appeared in a joint paper \cite{GustavssonJPeetreJ1977} of Jan Gustavsson
and Jaak. We intend to mention some more details about this method
and the way it is used in \cite{aaaCwKa} and beyond in a forthcoming
second more expanded version of \cite{CwikelMRochbergR2014Preprint}. 

As already mentioned, Nigel believed that ultimately a counterexample
will show that $\left(X_{0},X_{1}\right)\blacktriangleright\left(Y_{0},Y_{1}\right)$
does not hold in general. In parallel with his gentle and modest ways,
he was very ambitious, maybe even fiercely ambitious, and hoped very
much for the discovery of that counterexample. One of his wishes was
that it would show that the results of \cite{aaaCwKa} were just about
as strong as one could hope for. Well, in the almost 20 years since
the appearance of \cite{aaaCwKa}, as far as we know, only a small
number of additional partial answers to Question CIC have emerged.
Some of them will be described below. One could interpret this as
hinting that perhaps Nigel's wish may yet be granted, more or less. 

It is quite difficult by now to recall very much of the conversations
and ideas that were in the air during the writing of \cite{aaaCwKa},
beyond those which appear explicitly in the paper. A pity, because
perhaps they could yet be useful. Nigel liked to refer to some steps
of some arguments as ``the method of gliding humps''. Could that
mental picture be refined to guide us further?

It certainly seems worth going back to look again at many details
of \cite{aaaCwKa}. This is another one of the things we plan to do
in the future expanded version of \cite{CwikelMRochbergR2014Preprint}.
There is considerably more to be commented upon than could reasonably
fit into the limited framework of this document.

\section{Some developments since then}

\subsection{Some subsequently discovered additional cases where $\left(X_{0},X_{1}\right)\blacktriangleright\left(Y_{0},Y_{1}\right)$
holds. }

Here, as far as we know, are all the examples of couples $\left(X_{0},X_{1}\right)$
and $\left(Y_{0},Y_{1}\right)$ which have been shown to satisfy $\left(X_{0},X_{1}\right)\blacktriangleright\left(Y_{0},Y_{1}\right)$
since the appearance of \cite{aaaCwKa}.

In 2007 it was shown in \cite{CwikelJanson} that $\left(\ast,\ast\right)\blacktriangleright\left(FL^{\infty},FL_{1}^{\infty}\right)$.
Perhaps the main interest of this result lies in the hope that it
could maybe be some kind of first step towards showing that (\ref{eq:feid})
holds, and thus resolving Question CIC completely. 

It was particularly pleasing in that paper (and also in \cite{CwikelMJansonS2006})
to see Svante Janson returning to this field, when we recall that
his impressive work from some decades earlier (notably in the above-mentioned
paper\cite{JansonS1981}) is exceptionally important in interpolation
theory. We are also grateful to him for several very helpful conversations
in the intervening years.

In 2010 it was shown \cite{CwikelIntoLattices} that $\left(\ast,\ast\right)\blacktriangleright\left(Y_{0},Y_{1}\right)$
holds whenever $\left(Y_{0},Y_{1}\right)$ is a lattice couple for
which the underlying measure space is $\sigma$-finite and for which
some other mild condition is satisfied. This result was obtained,
roughly speaking, by taking a suitable variant of an adjoint operator
and applying a suitable variant of Schauder's classical theorem about
compact operators to the above-mentioned result (\ref{eq:felc}) of
\cite{aaaCwKa} (for which (as we intend to explain in more detail
in the future expanded version of \cite{CwikelMRochbergR2014Preprint})
there are no extra conditions required on the lattice couple). 

In the light of this result and of Theorem 11 of \cite[pp.~274--275]{aaaCwKa},
one is tempted to think that taking adjoints and using Schauder's
theorem provide an immediate and obvious way of obtaining a few more
partial new results. I.e., it seems reasonable to conjecture that
$\left(X_{0},X_{1}\right)\blacktriangleright\left(Y_{0},Y_{1}\right)$
if and only if their dual couples satisfy $\left(Y_{0}^{\prime},Y_{1}^{\prime}\right)\blacktriangleright\left(X_{0}^{\prime},X_{1}^{\prime}\right)$.
At this stage we only know how to show half of this conjecture, namely
that 
\begin{equation}
\left(Y_{0}^{\prime},Y_{1}^{\prime}\right)\blacktriangleright\left(X_{0}^{\prime},X_{1}^{\prime}\right)\mbox{ implies that }\left(X_{0},X_{1}\right)\blacktriangleright\left(Y_{0},Y_{1}\right)\label{eq:mbkztop}
\end{equation}
for any two regular couples $\left(X_{0},X_{1}\right)$ and $\left(Y_{0},Y_{1}\right)$.
(See Section 2 of the current preliminary version of \cite{CwikelMRochbergR2014Preprint}
for details.) This general fact seems to have been somehow overlooked
till now, though it was shown in one special setting in the proof
of Theorem 11 of \cite[pp.~274--275]{aaaCwKa}. The reverse implication
of (\ref{eq:mbkztop}) eludes us so far because, although all bounded
linear operators have adjoints, they do not always have ``pre-adjoints''.
In fact, whatever difficulty exists here is no smaller than the difficulty
of completely answering Question CIC. Why? Recalling the definitions
of the spaces $FL^{1}$, $FL_{1}^{1}$, $FC$ and $FC_{1}$ above,
a few lines before (\ref{eq:feid}), let us note that Calder\'on's
original result is sufficient to establish that $\left(\ell^{1}\left(FL^{1}\right),\ell^{1}\left(FL_{1}^{1}\right)\blacktriangleright\left(c_{0}\left(FC\right),c_{0}\left(FC_{1}\right)\right)\right)$.
If the reverse implication of (\ref{eq:mbkztop}) holds, then we can
almost immediately deduce that (\ref{eq:feid}) holds, which is equivalent
to answering ``yes'' to Question CIC. 

One immediate consequence of (\ref{eq:mbkztop}) is a new ``sister''
result for Theorem 9 of \cite[p.~271]{aaaCwKa}. That theorem told
us that $\left(X_{0},X_{1}\right)\blacktriangleright\left(Y_{0},Y_{1}\right)$
holds whenever $X_{0}$ is a $UMD$-space. Now we know that this holds
also whenever $Y_{0}$ is a $UMD$-space. (Again we refer to Section
2 of the current version of \cite{CwikelMRochbergR2014Preprint} for
more details.)

Some as yet unpublished results show that $\left(\ast,\ast\right)\blacktriangleright\left(E_{n,0},E_{n,1}\right)$
for every $n\in\mathbb{N}$ for some special particular couples $\left(E_{n,0},E_{n,1}\right)$,
which form a sequence which ``converges'' in some rather weak sense
to the couple which we would so like to have in place of them, namely
$\left(\ell^{\infty}\left(FL^{\infty}\right),\ell^{\infty}\left(FL_{1}^{\infty}\right)\right)$.
The elements in dense subclasses of all the spaces in all of these
couples can be conveniently visualized as sequences of holomorphic
functions of a complex variable on the ``unit annulus'' $\left\{ z\in\mathbb{C}:1\le\left|z\right|\le1\right\} $.
Complex interpolation requires us to consider holomorphic functions
of these functions. Thus we are dealing with sequences of holomorphic
functions of two complex variables.

\subsection{Other ``post-\cite{aaaCwKa}'' results related to Question CIC}

Soon after the appearance of \cite{aaaCwKa}, since Question CIC had
been defying attempts to answer it for so long, one of us, in collaboration
with Natan Kruglyak and Mieczys{\l}aw Mast{\l}o found it to be reasonable
to publish a ``toolbox'' paper \cite{CwikelMKrugljakNMastyloM1996}
identifying various special cases of Question CIC, whose resolution
would suffice to give a complete answer. An informal sequel to that
paper was also posted on the internet (\cite{CwikelMKrugljakNMastyloM9999InformaLetter}).
Subsequently, a remarkable calculation by Fedja Nazarov (in more formal
settings we should address him as Fedor) the details of which are
available at \cite{FejasNewExample}, hinted that perhaps there might
be a counterexample to a question which is closely related to Question
CIC, which is formulated as Question 2 on p.~362 of \cite{CwikelMKrugljakNMastyloM1996}.
(See also \cite{CwikelMKrugljakNMastyloM9999InformaLetter} or \cite[p.~168]{CwikelMJansonS2006}.)
The result of Fedja's calculation could also be interpreted as a first
vague hint that, even if the answer to Question CIC is affirmative,
there perhaps cannot be a generally valid quantitative version, in
terms of entropy numbers or covering numbers, of the property $\left(\ast,\ast\right)\blacktriangleright\left(\ast,\ast\right)$. 

It was remarked on p.~358 of \cite{CwikelMKrugljakNMastyloM1996}
that one can readily show that

\begin{equation}
\left(\ell^{\infty}\left(FL^{\infty}\right),\ell^{\infty}\left(FL_{1}^{\infty}\right)\blacktriangleright\left(\ell^{1}\left(FL^{1}\right),\ell^{1}\left(FL_{1}^{1}\right)\right)\right)\label{eq:tmi}
\end{equation}
via some results of Bill Johnson and Grothendieck. But (\ref{eq:tmi})
also follows from Theorem 11 of \cite[pp.~273--274]{aaaCwKa}. (Note
that (\ref{eq:tmi}) is a tantalizing ``mirror image'' of the result
(\ref{eq:feid}) that would completely answer Question CIC.) 

Some years later, some other authors also provided some additional
tools which might ultimately be helpful in further investigations
of Question CIC:

In 2000, Tomas Schonbek offered several interesting insights in his
paper \cite{SchonbekT2000}. He provided a more unified treatment
of results of \cite{aaaCwKa} and of other papers, clarifying the
central role which the property of equicontinuity plays in many of
them%
\footnote{Note, however, that at least one alternative strategy for trying to
answer Question CIC explicitly dispenses with using equicontinuity.
See Proposition 4 and the discussion preceding it on p.~359 of \cite{CwikelMKrugljakNMastyloM1996}.%
}. For example, his Theorem 3.5 on p.~1241 of \cite{SchonbekT2000}
gave an attractive alternative way of obtaining the result of Theorem
10 of \cite[p.~272]{aaaCwKa}. In Section 2 of his paper he indicated
a possible role for a variant of the Radon-Nikod\'ym Theorem. We
remark that, related to this, his Theorem 2.2 (p.~1234) might perhaps
turn out, in some contexts, to be able to replace or reinforce the
role played by Lemma 1 of \cite[p.~264]{aaaCwKa}.

In 2004, in their Theorem 3.2 of \cite[p.~71]{CobosFFernandezCabreraLMartinezA2004},
Fernando Cobos, Luz M.~Fern\'andez-Cabrera and Ant\'on Mart\'inez
gave necessary and sufficient conditions, somewhat in the style of
Lemma 5 of \cite[p.~268]{aaaCwKa} for a bounded operator $T:\left[A_{0},A_{1}\right]_{\theta}\to\left[B_{0},B_{1}\right]_{\theta}$
to be compact. (Their paper also dealt extensively with related issues
for other interpolation methods.) 

In 2006, Svante Janson and one of us \cite{CwikelMJansonS2006} briefly
reviewed some old and new possible strategies for dealing with Question
CIC. In particular, Section 3, pp.~166--168 of that paper contained
some observations about possible refinements, for this purpose, of
the method of ``infinite matrices of infinite matrices'' that was
used in \cite[Section 2, pp.~339--343]{CwikelM1992RCI} to prove
the ``extrapolation of compactness'' result mentioned here above,
near the end of Section \ref{sec:intro}. 

One year later Nigel, collaborating with Svitlana Mayboroda and Marius
Mitrea, obtained some interesting variants of that same earlier ``extrapolation
of compactness'' result. In Section 10 of \cite[pp.~160-163]{KaltonNMayborodaSMitreaM2007}
they ventured, as Nigel had so significantly done on so many other
occasions, beyond Banach spaces, into the much less ``comfortable''
setting of quasi-Banach spaces. They could deal with the cases where
these spaces are Besov or Triebel-Lizorkin spaces, and they applied
their results to partial differential equations.

Finally we mention a much more recent paper of J\"urgen Voigt \cite{VoigtJ2012}.
He considered the cases where the couples $\left(X_{0},X_{1}\right)$
and $\left(Y_{0},Y_{1}\right)$ satisfy either $X_{0}=X_{1}$ or $Y_{0}=Y_{1}$.
In such cases one immediately has $\left(X_{0},X_{1}\right)\blacktriangleright\left(Y_{0},Y_{1}\right)$
because of some classical results of Lions-Peetre. However the novelty
of the results in \cite{VoigtJ2012} is that, instead of dealing with
a single operator $T:X_{0}+X_{1}\to Y_{0}+Y_{1}$, they treat an operator
valued holomorphic function $T(z)$ defined for all $z$ in the strip
$\left\{ z\in\mathbb{C}:0\le\mathrm{Re}\, z\le1\right\} $ for which
$T(it):X_{0}\to Y_{0}$ is compact for all $t$ in some subset of
$\mathbb{R}$ having positive measure. (Of course operator valued
holomorphic functions, often also referred to as analytic families
of operators, have had many interactions with complex interpolation
theory, ever since they appeared in various special cases in \cite{HirschmanI1953}
and in \cite{SteinE1956}.)

\section{\label{sec:RelatedResults}Some related results, also for other interpolation
methods, and also for nonlinear operators}

It seems appropriate to say at least a few words about the large body
of results dealing with compact operators in the parallel contexts
of other interpolation methods. Fortunately we can refer to an extensive
and detailed survey \cite{CobosF2009} of many such results, written
by Fernando Cobos. (We are also grateful to Fernando for his many
other energetic initiatives to encourage research and interaction
in interpolation theory including but also well beyond issues of compactness,
via various conferences and meetings, and of course via his own numerous
results.) 

We recall that, as already mentioned above, the papers \cite{CobosKuhnSchonbek}
and \cite{CwikelM1992RCI} provide an affirmative answer for the analogue
of Question CIC for the Lions-Peetre real method of interpolation,
and that some special cases of this result had been known long before.
Once such qualitative compactness results have been obtained, it becomes
natural to seek quantitative refinements of them. Typically these
are expressed in terms of entropy numbers, or the measure of non-compactness.
Such results are amply described in \cite{CobosF2009}. Of course
there is a qualitative compactness result for the complex method too,
namely Calder\'on's seminal result, provided we remain within the
setting where the range couple $\left(Y_{0},Y_{1}\right)$ satisfies
Calder\'on's approximation condition. Very recently, Rados\l{}aw
Szwedek \cite{SzwedekR2015} has obtained quantitative results in
this setting. 

It is sometimes possible to obtain interpolation of compactness results
which apply simultaneously to a wide class of different interpolation
methods. Evgeniy Pustlynik \cite{PustylnikE2008} has done this when
the range couple $\left(Y_{0},Y_{1}\right)$ is a lattice couple having
certain extra properties. His results have some overlap with \cite{CwikelIntoLattices}.
Similar issues are addressed by Yuri Brudnyi in \cite{BrudnyiY2014}.

Finally let us briefly discuss interpolation of compact \textit{nonlinear}
operators by the real method. We confess to doing this mainly to give
us an opportunity to recall a special moment shared with Nigel. An
early paper \cite{CobosF1990} about this topic was written by Fernando
Cobos%
\footnote{But note that Fernando does not discuss nonlinear operators in \cite{CobosF2009}. %
}. It showed that, in some special cases, the affirmative answer to
the analogue of Question CIC for the real method of interpolation
(then already known from a preprint) could also be extended to the
case of nonlinear operators satisfying suitable properties, such as
Lipschitz conditions. We make no attempt to survey other research
on this topic, but only mention that when one of us met Nigel in Adelaide
in 2008 we casually started telling him about some joint work \cite{CwikelMIvtsanATadmorE9999}
extending the results of \cite{CobosF1990} and taking account of
some counterexamples we had found in \cite{cwikel-ivtsan}. Nigel,
being Nigel, instantaneously understood all factors in play, and immediately
guessed what condition we had needed to impose in order to obtain
our desired result.

\section{Anyway, why should we care about Question CIC?}

Of course we are attracted to Question CIC, as we are to other open
problems, in large part because of the fun and challenge of trying
to solve them, just as others are attracted to climbing Mount Everest.
However there seem to be some other good reasons for working on Question
CIC. 

Allowing ourselves to be a little ``carried away'', we can cautiously
claim that even the partial results obtained by Nigel and his coauthor
in \cite{aaaCwKa} might ultimately be relevant for confronting energy
crises. As evidence for this, we can consider the research by S\"oren
Bartels, Max Jensen and R\"udiger M\"uller, which is presented in
\cite{BartelsSJensenMMullerR2007} and then more formally in \cite{BartelsSJensenMMullerR2009}.
It deals with numerical solutions within a mathematical model relevant
for processes which extract oil from those underground reservoirs
where it is mixed with other fluids, and hints at the economic importance
of better understanding of such processes. Suddenly, in the middle
of both of these documents, there is an unexpected need for a tool
from interpolation theory, and Theorem 11 of \cite[p.~273--274]{aaaCwKa}
fulfills that need exactly. Nigel's coauthor saw fit to inform these
researchers that, in their particular setting, the seminal theorem
of Alberto Calder\'on, mentioned above (\cite{CalderonA1964} Sections
9.4 p.~118 and 29.4, pp.~137--138) would also suffice for their
purposes. One could also use results from the real interpolation method. 

Returning to the realm of ``pure'' or ``purer'' mathematics, we
can remark that there are references to \cite{aaaCwKa} in papers
dealing with such topics as compact Weyl operators, sub-Laplacians
of holomorphic $L^{p}$ type on exponential solvable groups, and boundary
value problems for Waldenfels operators. 

It could be argued that by now the answer to Question CIC is known
for all cases where the spaces $Y_{0}$ and $Y_{1}$ are ``reasonable''.
But mathematics and also its applications, sometime have a tendency
to move beyond what was previously considered ``reasonable''.

\bigskip{}

\textbf{\textit{Acknowledgment:}} We are grateful to Fernando Cobos,
Svante Janson, Mieczys\l{}aw Masty\l{}o and Mario Milman for some
very helpful correspondence.

\end{document}